\documentclass[12pt,reqno,a4wide]{amsart}

\usepackage{tikz} 
\usepackage{calrsfs}
\usepackage{mathrsfs}

\usepackage{array}
\usepackage{hyperref}

\usetikzlibrary{decorations.pathreplacing}
\usetikzlibrary{fit}

\usepackage{pgfplots}

\allowdisplaybreaks

\catcode`,\active

\catcode`\,12

\begin{document}
\bibliographystyle{plain}

%
%

	\title
	{Multi-edge trees and 3-coloured Motzkin paths: bijective issues}

	\author[H. Prodinger ]{Helmut Prodinger }
	\address{Department of Mathematics, University of Stellenbosch 7602, Stellenbosch, South Africa
	and
	Department of Mathematics and Mathematical Statistics,
	Umea University,
	907 36 Umea, 	 Sweden  }
	\email{hproding@sun.ac.za}

	
	\begin{abstract}
A bijection is given between multi-edge trees and 3-coloured Motzkin paths.\end{abstract}
	

\maketitle

\section{Introduction}

Multi-edge trees are like ordered (planar, plane, \dots) trees, but instead of edges there are multiple edges. When drawing such a tree,
instead of drawing, say 5 parallel edges, we just draw one edge and put the number 5 on it as a label. These trees were studied in
\cite{polish, HPW}. For the enumeration, one must count edges. The generating function $F(z)$ satisfies
\begin{equation*}
F(z)=\sum_{k\ge0}\Big(\frac{z}{1-z}F(z)\Big)^k=\frac1{1-\frac{z}{1-z}F(z)},
\end{equation*}
whence
\begin{equation*}
F(z)=\frac{1-z-\sqrt{1-6z+5z^2}}{2z}=1+z+3{z}^{2}+10{z}^{3}+36{z}^{4}+137{z}^{5}+543{z}^{6}+\cdots.
\end{equation*}
The coefficients form sequence A002212 in \cite{OEIS}.

A Motzkin path consists of up-steps, down-steps, and horizontal steps, see sequence A091965 in \cite{OEIS} and the references given there. As Dyck paths, they start at the origin and end, after $n$ steps again at the $x$-axis, but are not allowed to go below the $x$-axis. A 3-coloured Motzkin path is built as a Motzkin paths, but there
are 3 different types of horizontal steps, which we call \emph{red, green, blue}. The generating function $M(z)$ satisfies
\begin{equation*}
M(z)=1+3zM(z)+z^2M(z)^2=\frac{1-3z-\sqrt{1-6z+5z^2}}{2z^2}, \quad\text{or}\quad F(z)=1+zM(z).
\end{equation*}
So multi-edge trees with $N$ edges (counting the multiplicities) correspond to  3-coloured Motzkin paths of length $N-1$.

The purpose of this note is to describe a bijection. It transforms trees into paths, but all steps are reversible.

\section{The bijection}

As a first step, the multiplicities will be ignored, and the tree then has only $n$ edges. The standard translation of such tree into the world of Dyck paths,
which is in every book on combinatorics, leads to a Dyck path of length $2n$.
Then the Dyck path will transformed bijectively to a 2-coloured Motzkin path of length $n-1$ (the colours used are red and green).
This transformation plays a prominent role in \cite{Shapiro}, but is most likely much older. I believe that people like Viennot know this for 40 years.
I would be glad to get a proper historic account from the gentle readers.

The last step is then to use the third colour (blue) to deal with the multiplicities.

The first up-step and the last down-step of the Dyck path will be deleted. Then, the remaining $2n-2$ steps are coded pairwise into a 2-Motzkin path of length $n-1$:
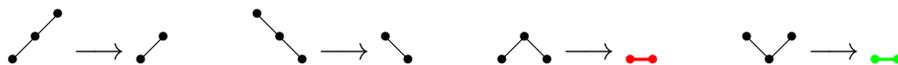
\begin{equation*}
\begin{tikzpicture}[scale=0.3]
					\path (0,0) node(x1) {\tiny$\bullet$} ;
\path (1,1) node(x2) {\tiny$\bullet$};
\path (2,2) node(x3) {\tiny$\bullet$};
	
	\draw (0,0) -- (2,2);

\end{tikzpicture}
\raisebox{0.5 em}{$\longrightarrow$}
\begin{tikzpicture}[scale=0.3]
	\path (0,0) node(x1) {\tiny$\bullet$} ;
	\path (1,1) node(x2) {\tiny$\bullet$};

	\draw (0,0) -- (1,1);

\end{tikzpicture}
\qquad
\begin{tikzpicture}[scale=0.3]
	\path (0,2) node(x1) {\tiny$\bullet$} ;
	\path (1,1) node(x2) {\tiny$\bullet$};
	\path (2,0) node(x3) {\tiny$\bullet$};
	
	\draw (0,2) -- (2,0);

\end{tikzpicture}
\raisebox{0.5 em}{$\longrightarrow$}
\begin{tikzpicture}[scale=0.3]
	\path (0,1) node(x1) {\tiny$\bullet$} ;
	\path (1,0) node(x2) {\tiny$\bullet$};
	
	\draw (0,1) -- (1,0);

\end{tikzpicture}
\qquad
\begin{tikzpicture}[scale=0.3]
	\path (0,0) node(x1) {\tiny$\bullet$} ;
	\path (1,1) node(x2) {\tiny$\bullet$};
	\path (2,0) node(x3) {\tiny$\bullet$};
	
	\draw (0,0) -- (1,1) -- (2,0);

\end{tikzpicture}
\raisebox{0.5 em}{$\longrightarrow$}
\begin{tikzpicture}[scale=0.3]
	\path (0,0) node[red](x1) {\tiny$\bullet$} ;
	\path (1,0) node[red](x2) {\tiny$\bullet$};
	
	\draw[red, very thick] (0,0) -- (1,0);

\end{tikzpicture}
\qquad
\begin{tikzpicture}[scale=0.3]
	\path (0,0) node(x1) {\tiny$\bullet$} ;
	\path (1,-1) node(x2) {\tiny$\bullet$};
	\path (2,0) node(x3) {\tiny$\bullet$};
	
	\draw (0,0) -- (1,-1) -- (2,0);

\end{tikzpicture}
\raisebox{0.5 em}{$\longrightarrow$}
\begin{tikzpicture}[scale=0.3]
	\path (0,0) node[green](x1) {\tiny$\bullet$} ;
	\path (1,0) node[green](x2) {\tiny$\bullet$};
	
	\draw[green, very thick] (0,0) -- (1,0);

\end{tikzpicture}
\end{equation*}

The last step is to deal with the multiplicities. If an edge is labelled with the number $a$, we will insert $a-1$ horizontal blue steps in the following way:
Since there are currently $n-1$ symbols in the path, we have $n$ possible positions to enter something (in the beginning, in the end, between symbols).
We go through the tree in pre-order, and enter the multiplicities one by one using the blue horizontal steps.

To make this procedure more clear, we prepared a list of 10 multi-edge trees with 3 edges, and the corresponding 3-Motzkin paths of length 2, with intermediate steps 
completely worked out:

\begin{center}
	\begin{table}[h]
		\begin{tabular}{c | c | c  |c}
			\text{Multi-edge tree  }&\text{Dyck path}&\text{2-Motzkin path}&\text{blue edges added}\\
			 	\hline\hline
				\begin{tikzpicture}[scale=0.5]
					\path (0,0) node(x1) {\tiny$\bullet$} ;
					\path (0,-1) node(x2) {\tiny$\bullet$};
					\path (0,-2) node(x3) {\tiny$\bullet$};
					\path (0,-3) node(x4) {\tiny$\bullet$};
					\draw (0,0) -- (0,-1)node[pos=0.5,left]{\tiny1} ;
					\draw (0,-1) -- (0,-2)node[pos=0.5,left]{\tiny1} ;
					\draw (0,-2) -- (0,-3)node[pos=0.5,left]{\tiny 1} ;

				\end{tikzpicture}
				 & \begin{tikzpicture}[scale=0.45]
				 	
				 	\draw (0,0) -- (3,3) --(6,0);

				 \end{tikzpicture}
				 & 
				 \begin{tikzpicture}[scale=0.45]
				 	
				 	\draw[thick] (0,0) -- (1,1) --(2,0);

				 \end{tikzpicture}
			 & \begin{tikzpicture}[scale=0.45]
			 	
			 	\draw[thick] (0,0) -- (1,1) --(2,0);

			 \end{tikzpicture}\\
			
			\hline

				\begin{tikzpicture}[scale=0.5]
					\path (0,0) node(x1) {\tiny$\bullet$} ;
					\path (0,-1) node(x2) {\tiny$\bullet$};
					\path (0,-2) node(x3) {\tiny$\bullet$};
				 
					\draw (0,0) -- (0,-1)node[pos=0.5,left]{\tiny2} ;
					\draw (0,-1) -- (0,-2)node[pos=0.5,left]{\tiny1} ;

				\end{tikzpicture}
			& \begin{tikzpicture}[scale=0.45]
				
				\draw (0,0) -- (2,2) --(4,0);

			\end{tikzpicture}
		& \begin{tikzpicture}[scale=0.45]
			
			\draw [red,thick](0,0) -- (1,0);

		\end{tikzpicture}
	& \begin{tikzpicture}[scale=0.45]
				\draw [blue,thick](0,0) -- (1,0);
		\draw [red,thick](1,0) -- (2,0);

	\end{tikzpicture}\\
			
			\hline
			\begin{tikzpicture}[scale=0.5]
				\path (0,0) node(x1) {\tiny$\bullet$} ;
				\path (0,-1) node(x2) {\tiny$\bullet$};
				\path (0,-2) node(x3) {\tiny$\bullet$};
				
				\draw (0,0) -- (0,-1)node[pos=0.5,left]{\tiny1} ;
				\draw (0,-1) -- (0,-2)node[pos=0.5,left]{\tiny2} ;

			\end{tikzpicture}
			& \begin{tikzpicture}[scale=0.45]
				
				\draw (0,0) -- (2,2) --(4,0);

			\end{tikzpicture}& \begin{tikzpicture}[scale=0.45]
			
			\draw [red,thick](0,0) -- (1,0);

		\end{tikzpicture}& \begin{tikzpicture}[scale=0.45]
		\draw [red,thick](0,0) -- (1,0);
		\draw [blue,thick](1,0) -- (2,0);

	\end{tikzpicture}\\
			
			\hline
			\begin{tikzpicture}[scale=0.5]
				\path (0,0) node(x1) {\tiny$\bullet$} ;
				\path (0,-1) node(x2) {\tiny$\bullet$};

				\draw (0,0) -- (0,-1)node[pos=0.5,left]{\tiny3} ;

			\end{tikzpicture}
			& \begin{tikzpicture}[scale=0.45]
				
				\draw (0,0) -- (1,1) --(2,0);

			\end{tikzpicture} & & \begin{tikzpicture}[scale=0.45]
			\draw [blue,thick](0,0) -- (2,0);

		\end{tikzpicture}\\
			
			\hline
			\begin{tikzpicture}[scale=0.5]
				\path (0,0) node(x1) {\tiny$\bullet$} ;
				\path (-1,-1) node(x2) {\tiny$\bullet$};
				\path (-2,-2) node(x3) {\tiny$\bullet$};
				\path (1,-1) node(x4) {\tiny$\bullet$};
				\draw (0,0) -- (-1,-1)node[pos=0.3,left]{\tiny1} ;
				\draw (-1,-1) -- (-2,-2)node[pos=0.3,left]{\tiny1} ;
				\draw (0,0) -- (1,-1)node[pos=0.3,right]{\tiny1} ;

			\end{tikzpicture}
			& \begin{tikzpicture}[scale=0.45]
				
				\draw (0,0) -- (2,2) --(4,0)--(5,1)--(6,0);

			\end{tikzpicture}& \begin{tikzpicture}[scale=0.45]
			
			\draw [red,thick](0,0) -- (1,0);
						\draw [green,thick](1,0) -- (2,0);
			
		\end{tikzpicture}& \begin{tikzpicture}[scale=0.45]
		
		\draw [red,thick](0,0) -- (1,0);
		\draw [green,thick](1,0) -- (2,0);
		
	\end{tikzpicture}\\
			
			\hline
			\begin{tikzpicture}[scale=0.5]
				\path (0,0) node(x1) {\tiny$\bullet$} ;
				\path (-1,-1) node(x2) {\tiny$\bullet$};
				
				\path (1,-1) node(x4) {\tiny$\bullet$};
				\draw (0,0) -- (-1,-1)node[pos=0.3,left]{\tiny2} ;
				
				\draw (0,0) -- (1,-1)node[pos=0.3,right]{\tiny1} ;

			\end{tikzpicture}
			& \begin{tikzpicture}[scale=0.45]
				
				\draw (0,0) -- (1,1) --(2,0)--(3,1)--(4,0);

			\end{tikzpicture}& \begin{tikzpicture}[scale=0.45]

			\draw [green,thick](0,0) -- (1,0);
			
		\end{tikzpicture}			& \begin{tikzpicture}[scale=0.45]
		
		\draw [blue,thick](0,0) -- (1,0);
		\draw [green,thick](1,0) -- (2,0);
		
	\end{tikzpicture}\\
			
			\hline
			\begin{tikzpicture}[scale=0.5]
				\path (0,0) node(x1) {\tiny$\bullet$} ;
				\path (-1,-1) node(x2) {\tiny$\bullet$};
				
				\path (1,-1) node(x4) {\tiny$\bullet$};
				\draw (0,0) -- (-1,-1)node[pos=0.3,left]{\tiny1} ;
				
				\draw (0,0) -- (1,-1)node[pos=0.3,right]{\tiny2} ;

			\end{tikzpicture}
			& \begin{tikzpicture}[scale=0.45]
				
				\draw (0,0) -- (1,1) --(2,0)--(3,1)--(4,0);

			\end{tikzpicture}& \begin{tikzpicture}[scale=0.45]

			\draw [green,thick](0,0) -- (1,0);
			
		\end{tikzpicture}& \begin{tikzpicture}[scale=0.45]
		
		\draw [green,thick](0,0) -- (1,0);
		\draw [blue,thick](1,0) -- (2,0);
		
	\end{tikzpicture}\\
			
			\hline
			\begin{tikzpicture}[scale=0.5]
				\path (0,0) node(x1) {\tiny$\bullet $} ;
				\path (-1,1) node(x2) {\tiny$\bullet$};
				\path (-2,2) node(x3) {\tiny$\bullet$};
				\path (-3,1) node(x4) {\tiny$\bullet$};
				\draw (0,0) -- (-1,1)node[pos=0.7,right]{\tiny1} ;
				\draw (-1,1) -- (-2,2)node[pos=0.7,right]{\tiny1} ;
				\draw (-2,2) -- (-3,1)node[pos=0.3,left]{\tiny1} ;

			\end{tikzpicture}
			& \begin{tikzpicture}[scale=0.45]
				
				\draw (0,0) -- (1,1) --(2,0)--(4,2)--(6,0);

			\end{tikzpicture}& \begin{tikzpicture}[scale=0.45]

			\draw [green,thick](0,0) -- (1,0);
\draw[red,thick](1,0)--(2,0);			
		\end{tikzpicture}& \begin{tikzpicture}[scale=0.45]

		\draw [green,thick](0,0) -- (1,0);
		\draw[red,thick](1,0)--(2,0);			
	\end{tikzpicture}\\
			
			\hline
			\begin{tikzpicture}[scale=0.5]
				\path (0,0) node(x1) {\tiny$\bullet$} ;
				\path (-1,-1) node(x2) {\tiny$\bullet$};
				\path (0,-1) node(x3) {\tiny$\bullet$};
				\path (1,-1) node(x4) {\tiny$\bullet$};
				\draw (0,0) -- (-1,-1)node[pos=0.3,left]{\tiny1} ;
				\draw (0,0) -- (0,-1)node[pos=0.6]{\tiny\;\;1} ;
				\draw (0,0) -- (1,-1)node[pos=0.3,right]{\tiny1} ;

			\end{tikzpicture}
			& \begin{tikzpicture}[scale=0.45]
				
				\draw (0,0) -- (1,1) --(2,0)--(3,1)--(4,0)--(5,1)--(6,0);

			\end{tikzpicture}& \begin{tikzpicture}[scale=0.45]

			\draw [green,thick](0,0) -- (1,0);
						\draw [green,thick](1,0) -- (2,0);
			
		\end{tikzpicture}& \begin{tikzpicture}[scale=0.45]

		\draw [green,thick](0,0) -- (1,0);
		\draw [green,thick](1,0) -- (2,0);
		
	\end{tikzpicture}\\
			
			\hline
			\begin{tikzpicture}[scale=0.5]
				\path (0,0) node(x1) {\tiny$\bullet$} ;
				\path (0,-1) node(x2) {\tiny$\bullet$};
				\path (-1,-2) node(x3) {\tiny$\bullet$};
				\path (1,-2) node(x4) {\tiny$\bullet$};
				\draw (0,0) -- (0,-1)node[pos=0.5,left]{\tiny1} ;
				\draw (0,-1) -- (1,-2)node[pos=0.5,right]{\tiny1} ;
				\draw (0,-1) -- (-1,-2)node[pos=0.5,left]{\tiny 1} ;

			\end{tikzpicture}
			& \begin{tikzpicture}[scale=0.45]
				
				\draw (0,0) -- (2,2) --(3,1)--(4,2)--(6,0);

			\end{tikzpicture}& \begin{tikzpicture}[scale=0.45]

			\draw [red,thick](0,0) -- (1,0);
			\draw [red,thick](1,0) -- (2,0);
			
		\end{tikzpicture}& \begin{tikzpicture}[scale=0.45]

		\draw [red,thick](0,0) -- (1,0);
		\draw [red,thick](1,0) -- (2,0);
		
	\end{tikzpicture}\\
			

		\end{tabular}
		
\caption{First row is a multi-edge tree with 3 edges, second row is the standard Dyck path (multiplicities ignored), third row is cutting off first and last step, and then translated pairs of steps, fourth row is inserting blued horizontal edges, according to multiplicities.}
	\end{table}	
\end{center}

\bibliographystyle{plain}


\end{document}